\documentclass{amsart} 
\usepackage{amssymb,latexsym, graphicx}

\newtheorem{theorem}{Theorem}[section]
\newtheorem{corollary}[theorem]{Corollary}

\theoremstyle{definition}
\newtheorem{definition}[theorem]{Definition}
\newtheorem{example}[theorem]{Example}

\theoremstyle{remark}
\newtheorem{remark}[theorem]{Remark}
\numberwithin{equation}{section}

\author{Anna Varvak}
\address{Soka University of America, Aliso Viejo, CA 92656}
\curraddr{Department of Mathematics,
      Brandeis University, Waltham, MA 02454}
\email{anka@brandeis.edu}

\title[Rook Numbers and the Normal Ordering Problem]
{Rook Numbers and the Normal Ordering Problem} 


\newcommand{\Sum}{\displaystyle\sum}
\newcommand{\Prod}{\displaystyle\prod}

\newcommand{\falling}[2]{{#1}^{\underline{#2}}}

\newcommand{\ms}
  {\vrule height 15pt
          depth  5pt
          width  0pt}
\newcommand{\df}{\ms\frac}

\newcommand{\ddx}{\frac{d}{dx}}
\newcommand{\ddt}{\frac{d}{dt}}

\DeclareMathOperator{\inv}{inv}

\DeclareMathOperator{\wt}{wt}
\DeclareMathOperator{\cross}{cross}
\DeclareMathOperator{\sep}{sep}

\def\vs#1{\vskip .#1 cm}

\catcode `!=11

\newdimen\squaresize 
\newdimen\thickness 
\newdimen\Thickness
\newdimen\ll! \newdimen \uu! \newdimen\dd! \newdimen \rr! \newdimen \temp!

\def\sq!#1#2#3#4#5{%
\ll!=#1 \uu!=#2 \dd!=#3 \rr!=#4
\setbox0=\hbox{%
 \temp!=\squaresize\advance\temp! by .5\uu!
 \rlap{\kern -.5\ll! 
 \vbox{\hrule height \temp! width#1 depth .5\dd!}}%
 \temp!=\squaresize\advance\temp! by -.5\uu!  
 \rlap{\raise\temp! 
 \vbox{\hrule height #2 width \squaresize}}%
 \rlap{\raise -.5\dd!
 \vbox{\hrule height #3 width \squaresize}}%
 \temp!=\squaresize\advance\temp! by .5\uu!
 \rlap{\kern \squaresize \kern-.5\rr! 
 \vbox{\hrule height \temp! width#4 depth .5\dd!}}%
 \rlap{\kern .5\squaresize\raise .5\squaresize
 \vbox to 0pt{\vss\hbox to 0pt{\hss $#5$\hss}\vss}}%
}
 \ht0=0pt \dp0=0pt \box0
}

\def\vsq!#1#2#3#4#5\endvsq!{\vbox to \squaresize{\hrule width\squaresize height 0pt%
\vss\sq!{#1}{#2}{#3}{#4}{#5}}}

\newdimen \LL! \newdimen \UU! \newdimen \DD! \newdimen \RR!

\def\vvsq!{\futurelet\next\vvvsq!}
\def\vvvsq!{\relax
  \ifx     \next l\LL!=\Thickness \let\continue=\skipnexttoken!
  \else\ifx\next u\UU!=\Thickness \let\continue=\skipnexttoken!
  \else\ifx\next d\DD!=\Thickness \let\continue=\skipnexttoken!
  \else\ifx\next r\RR!=\Thickness \let\continue=\skipnexttoken!
  \else\def\continue{\vsq!\LL!\UU!\DD!\RR!}%
  \fi\fi\fi\fi
  \continue}

\def\skipnexttoken!#1{\vvsq!}

\def\place#1#2#3{\vbox to 0pt{\vss
\rlap{\kern#1\squaresize
  \raise#2\squaresize\hbox{$#3$}}
\vss}}

\def\Young#1{\LL!=\thickness \UU!=\thickness \DD! = \thickness \RR! = \thickness
\vbox{\smallskip\offinterlineskip
\halign{&\vvsq! ## \endvsq!\cr #1}}}

\def\blank{\omit\hskip\squaresize}
\catcode `!=12

\begin{document}

\begin{abstract}
For an element $w$ in the Weyl algebra generated by $D$ and $U$ with relation $DU=UD+1$,
the normally ordered form is $w=\sum c_{i,j}U^iD^j$.  We demonstrate that the normal order 
coefficients $c_{i,j}$ of a word $w$ are rook numbers on a Ferrers board.  We use
this interpretation to give a new proof of the rook factorization theorem, which we use
to provide an explicit formula for the coefficients $c_{i,j}$.  We calculate the Weyl
binomial coefficients: normal order coefficients of the element $(D+U)^n$ in the Weyl
algebra. We extend these results to the $q$-analogue of the Weyl algebra.  We discuss further
generalizations using $i$-rook numbers.


\end{abstract}

\maketitle  


\section{Introduction}\label{S:Intro}

For an element $w$ in the Weyl algebra generated by $D$ and $U$ with relation $DU=UD+1$,
the normally ordered form is $w=\sum c_{i,j}U^iD^j$.  For example, in the algebra of
differential operators where $D=\ddx$ and $U$ acts as multiplication by $x$,
the operator $w$, applied to the polynomial $f(x)$, is expressed in the normally ordered
form as
$$
w\big( f(x) \big) = \Sum c_{i,j}\; x^i\, \frac{d^j f}{dx^j}(x) \;.
$$

The problem of finding explicit formulae for the normal order coefficients $c_{i,j}$
appears more frequently in the context where the Weyl algebra is the algebra of boson
operators~\cite{BPS, Katriel, Marburger, Mikh1, Mikh2, Schork}, generated by the creation and
annihilation operators typically denoted as $a^\dagger$ and $a$.
A boson is a type of particle like the light particle, the photon.  According to the theory
of quantum mechanics, the possible amount of energy that a particle can have is not
continuous but quantized, so there is the smallest amount of energy---the zeroth state,
frequently referred to as the ground state; there is the first state, which is the next
smallest amount of energy allowed, and so on.  The boson operators change the energy state
of the particle like the differential operators change the power of $x$.  While the
machanics of it are fascinating, for our purposes the assurance of the commutation relation 
$aa^\dagger - a^\dagger a = 1$ is sufficient.

As a hint of combinatorial interest in the problem of normal ordering, it has long been known
that the normal order coefficients of $(UD)^n$ are the Stirling numbers $S(n,k)$ of
the second kind.  These numbers can be defined algebraically by the formula
$$
x^n = \sum_{k=0}^{n} S(n,k) \; x(x-1)\cdots(x-k+1),
$$
and they also have a combinatorial interpretation of counting the number of ways
to partition a set of $n$ elements into $k$ subsets.

Mikhailov and Katriel \cite{Katriel, Mikh1, Mikh2} have extended the definition of the Stirling
numbers, finding explicit formulas for the normal ordered form of operators
such as $(D+U^r)^n$, and $(UD+U^r)^n$.
Recently, Blasiak, Penson, and Solomon~\cite{BPS} have generalized the Stirling numbers
even further to address the normal ordering problem for operators of the form $(U^rD^s)^n$.

The interpretation of the normal order coefficients of a word as rook numbers on a Ferrers
board was given by Navon~\cite{Navon} in 1973, but it requires the power of the rook
factorization theorem, presented two years later by Goldman, Joichi, and White~\cite{GJW} to
give an explicit formula.  Interestingly enough, the interpretation provides a proof of the
rook factorization theorem.

Here is an outline of the article.  In section~\ref{S:Basics}, we cover the basic
definitions concerning Ferrers boards and rook numbers.  We demonstrate that the
normal order coefficents of a word are rook numbers on a Ferrers board in
section~\ref{S:NOPWord}, and give an explicit formula for them in section~\ref{S:CompNOC},
together with a new proof of the rook factorization theorem.  We discuss the Weyl binomial
coefficients in section~\ref{S:WeylBinom}.  We extend the above results to their $q$-analogue in
section~\ref{S:q-analogue}.  Finally, we discuss some further generalizations in
section~\ref{S:i-analogue}.

\section{Definitions concerning rook numbers and Ferrers boards}\label{S:Basics}

For $n$ a positive integer, we denote by $[n]$ the set $\{ 1,2,\dots, n\}$. A
\emph{board} is a subset of $[n] \times [m]$, where $n$ and $m$ are positive
integers.  Intuitively, we think of a board as an array of squares arranged in rows
and columns.  An element $(i,j) \in B$ is then represented by a square in the 
$i$th column and $j$th row.  It will be convenient to consider columns numbered from left
to right, and rows numbered from top to bottom, so that the square $(1,1)$ appears in the
top left corner.

A board $B$ is a \emph{Ferrers} board if there is a non-increasing sequence of positive
integers
$h(B)=(h_1, h_2, \dots, h_n)$  such that 
$B = \{(i,j) \;|\; i \le n \mbox{ and } j\le h_i \}$. 
Intuitively, a Ferrers board is a
board made up of adjacent solid columns with a common upper edge, such that the heights
of the columns from left to right form a non-increasing sequence. 

\begin{example}\label{E:FerrersBoard}
A Ferrers board with height sequence $(4,4,3,1,1)$ can be visually represented as
$$
\Young{
       &       &       &       & \cr
       &       &       &\blank &\blank \cr
       &       &       &\blank &\blank \cr
       &       &\blank &\blank &\blank \cr
}
$$
\end{example}

The connection between Ferrers boards and words composed of two letters is as follows. 
We represent the letter $D$ as a step to the right, and the letter $U$ as a step up.  
The resulting path outlines a Ferrers board.

\begin{example}
The word $w=DDUDUUDDU$ outlines the Ferrers board in example~\ref{E:FerrersBoard}.
This is easy to see from the path representing $w$:

$$
\Young{
\blank &\blank &\blank &lur    &lu \cr
\blank &\blank &lud    &\blank &\blank \cr
\blank &\blank &lu     &\blank &\blank \cr
lur    &lu     &\blank &\blank &\blank \cr
}
$$

Note that a word $w'=U^i w D^j$ outlines the same Ferrers board for any
nonnegative integers $i$ and $j$.  If the Ferrers board $B$ is contained in a rectangle with
$n$ columns and $m$ rows, then there is a unique word with $n$ $D$'s and $m$ $U$'s that
outlines
$B$.
\end{example}

We denote the Ferrers board outlined by the word $w$ by $B_w$.

For a board $B$, let $r_k(B)$ denote the number of ways of marking $k$ squares of the
board $B$, no two in the same row or column.  
In chess terminology, we are placing $k$ rook pieces on the board $B$ in non-attacking
positions. The number $r_k(B)$ is called the $k$th \emph{rook number} of $B$.

\section{Normal order coefficients of a word}\label{S:NOPWord}

Recall that the Weyl algebra is the algebra generated by $D$ and $U$,
with the commutation relation $DU= UD+1$.  

\begin{definition}\label{D:NormalOrderedForm}
For $w$ an element in the Weyl algebra, the \emph{normally ordered form} of $w$ is
the sum
$$
w = \Sum_{i,j} c_{i,j} U^i D^j,
$$
where in each term the $D$'s appear to the right of the $U$'s.  The numbers
$c_{i,j}$ are the \emph{normal order coefficients} of $w$.
\end{definition}

We call $w$ a \emph{word} if $w$ has a representation $w=w_1w_2 \dots w_n$,
where $w_i \in \{D,U\}$.  We demonstrate that the normal order coefficients of a word $w$
are rook numbers on the Ferrers board outlined by $w$.  This combinatorial interpretation
was originally given by Navon~\cite{Navon}.

\begin{theorem}\label{T:main} {\bf Normally Ordered Word}

Let the element $w$ in the Weyl algebra be a word composed of $n$ $D$'s and $m$ $U$'s.
Then
$$
w = \sum_{k=0}^{n} r_k(B_w) U^{m-k}D^{n-k}.
$$
\end{theorem}

\begin{proof}
It is easy to see that the terms in the normally ordered form of $w$ are
$U^{m-k}D^{n-k}$, where $k=0,1,\dots, \min(m,n)$.  Every time we replace $DU$
with $UD+1$ and expand the result, one of the new terms retains the same number of $D$'s
and $U$'s as before, and the other term has one fewer of each.  

By definition, the normal order coefficient $c_{m-k, n-k}$ is the number of terms
$U^{m-k}D^{n-k}$ in the normally ordered form of $w$, which is obtained by successively
replacing $DU$ with $UD+1$ and expanding.  For the sake of consistency, we
always choose to replace the rightmost $DU$.

We can regard the terms as words.  Then the normal order coefficient $c_{m-k, n-k}$ is
the number of ways to get the word $U^{m-k}D^{n-k}$ from the word $w$, by successively
replacing the rightmost $DU$ with either $UD$ or $1$ (that is, deleting it), choosing
to do the latter $k$ times.

We now can consider the substitutions in terms of the outlined Ferrers boards.  The
rightmost $DU$ outlines the rightmost inner corner square of the board.  Replacing
$DU$ with $UD$ amounts to deleting that square, whereas deleting the $DU$ amounts to
deleting the square together with its row and column.  Therefore the normal order
coefficient $c_{m-k, n-k}$ is the number of ways to reduce the Ferrers board $B_w$
outlined by $w$ to the trivial board by successively  deleting the rightmost inner corner
square either alone, or together with its row and column, choosing to do the latter $k$
times.

The $k$ squares that are deleted together with their rows and columns cannot share either
a row or a column.  So the normal order coefficient $c_{m-k, n-k}$ is the number
of ways to mark $k$ squares on the Ferrers board $B_w$ outlined by the word $w$, no two
in the same row or column.  This is exactly the $k$th rook number $r_k(B_w)$.

Finally, since $r_k(B_w)=0$ for $k>\min(m,n)$, we let the sum range
from~$1$ to~$n$.
\end{proof}

\begin{remark}

As mentioned in the introduction, the Stirling numbers $S(n,k)$ of the second kind
are the normal order coefficients of the word $(UD)^n$.
Mikhailov \cite{Mikh1} defined, in a purely algebraic way, a more
generalized version of the Stirling numbers to find the normal ordered form of
operators of the form $(U^r+D)^n$.  In a recent paper unrelated to the normal ordering
problem, Lang \cite{Lang} studied a similar generalization of the Stirling numbers,
finding combinatorial interpretations for certain particular cases.  Recently
Blasiak, Penson, and Solomon~\cite{BPS} introduced the generalized Stirling
numbers of the second kind, denoted $S_{r,s}(n,k)$ for $r \ge s \ge 0$, defined by the
relation
$$
\left( U^r \, D^s \right)^n 
  = U^{n(r-s)} \sum_{k=s}^{ns} S_{r,s}(n,k) U^k \, D^k.
$$
The standard Stirling numbers of the second kind are $S_{1,1}(n,k)$, and the generalized
Stirling numbers of Mikhailov are $S_{r,1}(n,k)$.

We define the \emph{staircase board} $J_{r,s,n}$ to be the Ferrers board outlined by the word
$(U^rD^s)^n$.

\begin{corollary}
$$
S_{r,s}(n,k) = r_{ns-k}\left(J_{r,s,n} \right).
$$
\end{corollary}

\end{remark}

\begin{remark}
We can easily adapt the proof of Theorem~\ref{T:main} to work in the case where
the algebra generated by $D$ and $U$ has the commutation relation $DU=UD+c$.  For a word $w$,
we get the normal ordering form of $w$ by successively replacing $DU$ by $UD+c$, and
expanding.  Just as before, we consider $w$ as a word in the letters $D$, $U$, and the
substitution as a choice of either replacing the rightmost $DU$ by $UD$ or deleting it, but
the choice of deleting is weighted by $c$. In terms of the associated Ferrers board, we weight
each placement of a rook by $c$. So we get   
$$
w = \sum_{k=0}^{n} c^k r_k(B_w) U^{m-k}D^{n-k}.
$$
We should note that this algebra is isomorphic to the Weyl algebra, because if
$DU-UD=1$, then $D(cU)-(cU)D=c$.

We can also assign a weight to the choice of replacing $DU$ by $UD$, and thus
extend the result to algebra with the relation $DU=qUD+1$.  The algebra with this
relation is know as the $q$-Weyl algebra, and is of interest
both to combinatorialists and to physicists.  To the latter, because such algebras are
models for $q$-degenerate bosonic operators~\cite{Schork}.  To the former, because it
involves the
$q$-analogue of rook numbers~\cite{GR}.  We discuss this case in detail in
section~\ref{S:q-analogue}.

First, we show how Theorem~\ref{T:main}
allows for a new proof of the Rook Factorization Theorem~\cite{GJW}, which in turn
leads to an explicit formula for computing the normal order coefficients of a word.
\end{remark}

\section{Computing the normal order coefficients}\label{S:CompNOC}

For a general board $B$, rook numbers can be computed recursively~\cite{Riordan}.  Choose
a square of
$B$, and let $B_1$ be the board obtained from $B$ by deleting that square, and let $B_2$
be the board obtained from $B$ by deleting the square together with its row and column.
Then  $r_k(B) = r_k(B_1)+r_{k-1}(B_2)$, reflecting the fact we may or may not
mark the square in question.

There are better methods for calculating rook numbers on Ferrers boards, owing to the fact
that the generating function of rook numbers on a Ferrers board, expressed in terms of
falling factorials, completely factors. 

We define the $k$th \emph{falling factorial} of $x$ by
$$\falling{x}{k}=x(x-1)\cdots(x-k+1).$$
Goldman et al.~\cite{GJW} show that
the factorial rook polynomial 
$\sum_{k=0}^{n}r_k(B) \, \falling{x}{n-k}\;$
of a Ferrers board  is a product of linear factors.

\begin{theorem}\label{T:RFT1}{\bf Rook Factorization Theorem}

For a Ferrers board $B$ with column heights $h(B) = (h_1, \dots, h_n)$,
$$
\sum_{k=0}^{n}r_{k}(B) \falling{x}{n-k} = \prod_{i=1}^{n}(x+h_i-n+i)
$$
\end{theorem}

We provide a new proof the Rook Factorization Theorem, using Theorem~\ref{T:main}.

\begin{proof}
Let $w$ be the word with $n$ $D$'s and $h_1$ $U$'s that outlines the Ferrers board~$B$.
By Theorem~\ref{T:main}, 
$$
w = \sum_{k=0}^{n} r_k(B_w) U^{h_1-k}D^{n-k}.
$$
as an element in the Weyl algebra.

We consider a particular manifestation of the Weyl algebra as the algebra of operators
generated by $D=\frac{d}{dt}$ and
$U$= multiplication by $t$, acting on functions in the variable $t$.  So
$$
w = \sum_{k=0}^{n} r_k(B_w) t^{h_1-k}\left( \frac{d}{dt} \right)^{n-k}.
$$
We apply both sides of the equation to $t^x$, where $x$ is a real number.

Since $\big( \frac{d}{dt} \big)^{n-k}(t^x) = x(x-1) \cdots (x-(n-k)+1)
t^{x-n+k}=\falling{x}{n-k} \; t^{x-n+k}$, the right hand side is
$$
\sum_{k=0}^{n} r_k(B) \, \falling{x}{n-k} \; t^{x-n+h_1}.
$$
On the left-hand side we get the product of the following factors.  
The $j$th application of $D$ gives the factor of
$x+a_U-a_D$, where
$a_U$ the number of times $U$ was previously applied, and $a_D$ the number of times $D$
was previously applied.  There are $j-1$ $D$'s to the right of the $j$th $D$, so
$a_D=j-1$. The $j$th $D$ from the right is the $(n-j+1)$st $D$ from the left, so
$a_U=h_{n-j+1}$. Therefore the left-hand side is 
$$
t^{x-n+h_1}\prod_{j=1}^{n}(x+h_{n-j+1}-j+1). 
$$
If we let $i=n-j+1$, then the left-hand side is
$$
t^{x-n+h_1}\prod_{i=1}^{n}(x+h_i-n+i). 
$$
Now we set $t=1$ to get the desired result.
\end{proof}

\begin{example}
For $w = DDUUUDDUD$, by Theorem~\ref{T:main}, 
\begin{align*}
\ddt \,\ddt\,t \cdot t \cdot t \cdot \ddt \,\ddt\, t\cdot \ddt \;\big(t^x\big) 
  &= \Sum_{k=0}^{5} r_k(B_w)\;t^{4-k}\,\left( \ddt \right)^{5-k}\;\big(t^x\big)\\
x\cdot (x+1)\cdot (x-1)\cdot x\cdot x\cdot \big(t^{x-1}\big) 
  &= \Sum_{k=0}^{5} r_k(B_w)\;\falling{x}{5-k}\cdot \big(t^{x-1}\big)\\
x\cdot (x+1)\cdot (x-1)\cdot x\cdot x  
  &= \Sum_{k=0}^{5} r_k(B_w)\;\falling{x}{5-k} \;.
\end{align*}
The left hand side is the complete factorization of the factorial rook polynomial of $B_w$.

\end{example}

The falling factorials $1, x, x(x-1), \dots$ form a basis of polynomials in $x$.  If
$P(x) = \sum_{k=0}^{n} p_k \, \falling{x}{k},$
it is well known that the coefficients are
$p_k = \df{1}{k!} \Delta^{k} P(x)\Big|_{x=0},$
where $\Delta$ is the difference operator defined by $\Delta P(x) = P(x+1)-P(x)$.  
Explicitly~\cite{StVol1}, the coefficients are
$$
p_k = \df{1}{k!} \Sum_{i=0}^{k} (-1)^{k-i}{k \choose i}\,P(i)\;.
$$

\begin{corollary}{\bf Computing the Normal Order Coefficients of a Word}\label{C:comp}

For $w$ a word in the Weyl algebra composed of $n$ $D$'s and $m$ $U$'s,
let $P(x) = \prod_{i=1}^{n}(x+h_i-n+i)$, where $h_1, h_2, \dots, h_n$ are the column
heights of the Ferrers board outlined by $w$.  Then
$$
w = \sum_{k=0}^{n} r_k U^{m-k}D^{n-k},
$$
where
$$
r_k = \df{1}{(n-k)!} \Sum_{i=0}^{n-k} (-1)^{n-k-i}{n-k \choose i}\,P(i)\;.
$$
\end{corollary}

\section{Weyl binomial coefficients}\label{S:WeylBinom}

The Weyl binomial coefficient ${n\choose m}_k$ is the coefficient of the term
$U^{n-m-k}D^{m-k}$ in $(D+U)^n$, where the commutation relation is $DU=UD+1$.  The binomial
product
$(D+U)^n$ is the sum of all words in letters $D$ and $U$ of length $n$, and the normally
ordered term
$U^{n-m-k}D^{m-k}$ comes from all words with $n-m$ $U$'s and $m$ $D$'s, where $k$ pairs of
$D$ and $U$ are deleted during the normal ordering.  Each of these words outlines a unique
Ferrers board with at most $m$ columns of height at most $(n-m)$. Therefore the Weyl
binomial coefficient
${n\choose m}_k$ can be expressed as the sum of the $k$th rook numbers over all Ferrers
boards contained in the $m$-by-$(n-m)$ rectangle.

The classical binomial coefficient ${n \choose m}$, and its $q$-analogue ${n\brack m}$,
can be similarly expressed in terms of Ferrers boards.
The binomial coefficient ${n \choose m}$ is the coefficient of the term $U^{n-m}D^{m}$ in
$(D+U)^n$, where the commutation relation is $DU=UD$.
The term $U^{n-m}D^{m}$ is the normally ordered form of any word $w$ with $n-m$ $U$'s and $m$
$D$'s.  Since all the letters commute, there is only one way to normally
order a word, so the binomial coefficient ${n \choose m}$ counts the words with $n-m$ 
$U$'s and $m$ $D$'s.  Since there is a bijection between the set of such words and the set
of Ferrers boards contained in $m$-by-$(n-m)$ rectangle, ${n \choose m}$ counts such Ferrers
boards.  Therefore
$$
{n \choose m} = \sum_{B \subseteq [m]\times[n-m]} 1.
$$

Similarly, the $q$-binomial coefficient ${n\brack m}$ is the coefficient of the term
$U^{n-m}D^{m}$ in $(D+U)^n$, where the commutation relation is $DU=qUD$.  Again, the term
$U^{n-m}D^{m}$ comes from any word $w$ with $n-m$ $U$'s and $m$ $D$'s, which outlines the
Ferrers board $B_w$ contained in the rectangle $[m]\times[n-m]$.  The relation $DU=qUD$
specifies that each square of
$B_w$ has weight $q$, so  $w = q^{|B_w|}U^{n-m}D^{m}$.  In other words,
$$
{n \brack m} = \sum_{B \subseteq [m]\times[n-m]} q^{|B|}.
$$

\begin{theorem}\label{T:WeylBinCoef}
Let $k \le m$ be an integer.  Then  
$$
\sum_{B \subseteq [m]\times[n-m]} r_k(B)
={n\choose m}_k = \frac{n!}{ 2^k\,k!\,(m-k)!\,(n-m-k)!}
$$
\end{theorem}

\begin{proof}
Any Ferrers board $B$ in $[m]\times[n-m]$ is $B_w$ for a particular
word
$w$ with $n-m$ $U$'s and $m$ $D$'s.  From the proof of Theorem~\ref{T:main}, the number
$r_k(B_w)$ is the number of ways to get the word $U^{n-m-k}D^{m-k}$ from the word $w$, by
successively either replacing $DU$ with $UD$ or
deleting it, choosing to do the latter $k$ times.
Equivalently, $r_k(B_w)$ is the number of ways to mark $k$ pairs of
the letters $D$ and $U$ in the word $w$, such that in each pair the $D$ appears to the
left of the $U$. The marked pairs are deleted, and the rest of the letters commute into
the normally ordered form.

We therefore count the number of ways to construct words with $n-m$ $U$'s and $m$ $D$'s,
with $k$ marked pairs of the letters $D$ and $U$, the former to the left of the latter.
We begin with $n$ spaces for the letters in $w$.  Choose $n-m-k$ of these to be $U$,
and $m-k$ to be $D$.  There are
$$
{n \choose n-m-k, \; m-k} = \frac{n!}{(n-m-k)!\,(m-k)!\,(2k)!}
$$
ways to do so.  For the $2k$ remaining spaces, we pair them, forming $k$ pairs.
There are $(2k-1) \cdot (2k-3) \cdots 5\cdot 3\cdot 1$ ways to do this.
For each pair, let the space on the left be $D$, and the space on the right be $U$.

By this construction, the sum of the $k$th rook numbers over all Ferrers boards contained
in the $m$-by-$(n-m)$ rectangle is
$$
\frac{(2k-1) \cdot (2k-3) \cdots 5\cdot 3\cdot 1}{(2k)!}\frac{n!}{(n-m-k)!\,(m-k)!}
= \frac{n!}{ 2^k\,k!\,(m-k)!\,(n-m-k)!}.
$$

\end{proof}

Since $(D+U)^n$ is the sum of all words composed of letters $D$ and $U$ of length $n$,
we have a formula for normal ordering of $(D+U)^n$, as shown by Mikhailov in
\cite{Mikh1}.

\begin{corollary}\label{C:binom}
$$
(D+U)^n = 
\sum_{m=0}^{n} 
\sum_{k=0}^{\min(m, n-m)} \frac{n!}{ 2^k\,k!\,(m-k)!\,(n-m-k)!} U^{m-k}D^{n-m-k}
$$
\end{corollary}

\begin{remark}
The Weyl binomial coefficients obey the recursive formula
$$
{n \choose m}_k = {n-1 \choose m}_k + {n-1 \choose m-1}_k + m {n-2 \choose m-1}_{k-1} \;,
$$
with boundary conditions ${1 \choose 0}_0={1\choose 1}_0=1$, ${n\choose m}_{k-1}=0$.
Consider the pairs $(B,C)$, where $B$ is a Ferrers board in $[m]\times[n-m]$ and
$C$ is a placement of $k$ rooks on $B$.  The Weyl binomial coefficient 
counts the number of such pairs. 
The set of such pairs is a disjoint union of three sets: one where the height $h_1$ of the
first column $B$ is strictly less than ${n-m}$, one where $h_1=m$ and $C$ doesn't place a
rook in the first column, and one where $h_1=m$ and $C$ places a rook in the first column.
The recursive formula follows. 

The first two terms in the recursive formula are the same as for the classical binomial
coefficients.  In fact, the Weyl binomial coefficients can be expressed in terms of
classical coefficients as follows.  

\begin{corollary}
Let $C(y) = \sum_{k\ge 0} {n\choose k}\, {y^k\over k!}$ be the exponential
generating function of the binomial coefficients.  Then the ordinary generating function of
the Weyl binomial coefficients is
$$
\Sum_{k\ge 0} {n\choose m}_k\,x^k 
=\left. \left(\df{d}{dy}\right)^{n-m}\,C(y)\;\right|_{y=\df{x}{2}} \;.
$$
\end{corollary}

\begin{proof}
\begin{align*}
C(y) &= \Sum_{k\ge 0} {n\choose k}\, {y^k\over k!} \\
     &= \Sum_{k\ge 0} {n\choose n-k}\, {y^k\over k!} \;, \\
\end{align*}
so
$$
\left(\df{d}{dy}\right)^{n-m}\,C(y) 
      = \Sum_{k\ge 0} {n\choose m-k}\, {y^k\over k!}\;, \\
$$
therefore
\begin{align*}
\left. \left(\df{d}{dy}\right)^{n-m}\,C(y)\;\right|_{y=\df{x}{2}}
               &= \Sum_{k\ge 0} {n\choose m-k}\,{1\over k!} {x^k\over 2^k} \\
               &= \Sum_{k\ge 0} \frac{n!}{ 2^k\,k!\,(m-k)!\,(n-m-k)!} \, x^k \;.
\end{align*}
\end{proof}

\end{remark}

\section{The $q$-analogue}\label{S:q-analogue}

We extend the combinatorial interpretation of the normal order coefficients to the $q$-Weyl
algebra: the algebra with two generators $D$ and $U$, and the relation $DU=qUD+1$.

The commutation relation twisted by $q$ comes up in physics as the relation obeyed 
by the creation and annihilation operators of $q$-deformed bosons~\cite{Schork}.  The problem
of normally ordering these operators has been studied by Katriel~\cite{Katriel}, and
recently by Schork~\cite{Schork}.

The basic idea of the $q$-analogue of numbers is that the polynomial
$q^0+q^1+q^2+\cdots+q^{n-1}$ plays the role of the positive integer $n$.  We denote the
$q$-analogue of $n$ by $[n]_q$. Since $1+q+q^2+\cdots+q^{n-1} = \frac{1-q^n}{1-q}$, we can 
extend the definition of the $q$-analogue to  all numbers $t$ by defining
$$
[t]_q := \frac{1-q^t}{1-q}.
$$
The $q$-analogue of the derivative $\frac{d}{dt}$ acting on the ring of polynomials
in $t$ is defined as
$$
D_q\,f(t) := \df{f(qt)-f(t)}{(q-1)t}\;,
$$
and it is easy to check that $D_q(t^n)=[n]_q\,t^{n-1}$.  For a good exposition of the 
$q$-analogue of the derivative, we refer the reader to ``Quantum Calculus" by Kac and 
Cheung~\cite{KC}.

If we let $D = D_q$, and let $U$ be the operator acting by
multiplication by $t$, then the algebra generated
by $D$ and $U$ has the relation $DU=qUD+1$, and is therefore the $q$-Weyl algebra.

Let the element $w$ in the $q$-Weyl algebra be a word composed of the letters $D$ and $U$. 
We adapt the proof of Theorem~\ref{T:main} to find the normal order
coefficients of $w$. 

In terms of algebraic operations, we get the normal ordering form of
$w$ by successively replacing $DU$ by $qUD+1$, and expanding.  As before, we can consider
$w$ as a formal word in the letters $D$ and $U$, and the substitution as a choice of either
replacing the rightmost $DU$ by $UD$, weighting this choice by $q$, or deleting the
rightmost $DU$. In terms of the Ferrers board $B_w$ outlined by $w$, we assign the weight
$q$ to each square that doesn't have a rook either on it, below it in the same column, or
to the right of it in the same row.  If we consider the weight of a rook placement to be
the product of the weights of all squares of the board, such weights of rook placements
describe exactly the $q$-rook numbers of Garsia and Remmel~\cite{GR}.

\begin{definition}
Let $B$ be a board, and denote by $C_k(B)$ the collection of all placements of $k$
marked squares (rooks) on $B$, no two in the same row or column.  We define the $k$th
$q$-rook number of $B$ to be
$$
R_k(B,q) = \sum_{C \in C_k(B)} q^{\inv(C)},
$$
where $\inv(C)$ is the number of squares in the placement $C$ 
that do not have a rook either on them, below them in the same column, or
to the right of them in the same row.
\end{definition}

\begin{remark}
To clarify the statistic $\inv(C)$, we demonstrate by an example.  Suppose $B$
is a Ferrers board with column heights $h(B)=(4,4,3,1)$, and we have the following
placement $C$ of two rooks
$$
\Young{
       &        &        &         \cr
       &        & \times &\blank   \cr
       & \times &        &\blank   \cr
       &        &\blank  &\blank   \cr
}
$$
If we mark with a dot the squares above or to the left of a rook, and with a circle
the rest of the squares, we get
$$
\Young{
\circ  & \cdot  & \cdot  &\circ    \cr
\cdot  & \cdot  & \times &\blank   \cr
\cdot  & \times & \circ  &\blank   \cr
\circ  & \circ  &\blank  &\blank   \cr
}
$$
Then $\inv(C)$ is the number of squares marked with a circle, which in this example is~5.

The statistic $\inv(C)$ is a generalization of the inversion statistic
on permutations.  Given a permutation $\sigma=(\sigma_1,\dots, \sigma_n)$, we get a
placement $C$ of $n$ rooks on an $n$-by-$n$ board where the rook in column $i$ is placed
in row $\sigma_i$.  Then each square marked with a circle has a rook to the left of
it and a rook above it, so the square corresponds to an inversion pair $i<j$ 
such that $\sigma_i > \sigma_j$.  So in this case, $\inv(C)$ is the number of inversions of
$\sigma$.  
\end{remark}

Analogous to Theorem~\ref{T:main}, the normal order coefficients of a word
$w$ are the $q$-rook numbers of the Ferrers board outlined by $w$.

\begin{theorem}\label{T:qmain}
Let the element $w$ in the $q$-Weyl algebra be a word composed of $n$ $D$'s and $m$ $U$'s.
Then
$$
w = \sum_{k=0}^{n} R_k(B_w, q) U^{m-k}D^{n-k}.
$$
\end{theorem}

The Factorization Theorem for $q$-rook
numbers~\cite{GR} can be proved as a corollary.

\begin{theorem} {\bf Factorization Theorem for $q$-Rook Numbers}

For a Ferrers board $B$ with column heights $h(B) = (h_1, \dots, h_n)$,
$$
\sum_{k=0}^{n}R_k(B,q) [x]_q\,[x-1]_q\,\cdots[x-(n-k)+1]_q\, 
=\prod_{i=1}^{n}[x+h_i-n+i]_q\,
$$
\end{theorem}

The proof is exactly the same as for
Theorem~\ref{T:RFT1}, replacing the real numbers, the falling factorials, and the
derivative with their $q$-analogue, and using the fact that $D_q(t^n)=[n]_q\,t^{n-1}$.


We provide a formula for the normal order coefficients of a word as in Corollary~\ref{C:comp},
but with the following alterations~\cite{GR}.  If 
$$
P(x) = \Sum_{k=0}^{n} p_k\,[x]_q\,[x-1]_q\,\dots[x-k+1]_q\;,
$$
then
$$
p_k = \df{1}{[k]_q\,!} \Delta_q^k\,P(x)\Big|_{x=0} \;,
$$
where 
$\Delta_q^k\,P(x)=\Prod_{i=1}^{k}\big( P(x+1)-q^{i-1}P(x) \big)$
and $[k]_q\,!=[k]_q\,[k-1]_q\,\dots[1]_q\,$.

\begin{corollary}\label{C:q-comp}

For $w$ a word in the $q$-Weyl algebra composed of $n$ $D$'s and $m$ $U$'s,
let $P(x) = \prod_{i=1}^{n}[x+h_i-n+i]_q\,$, where $h_1, h_2, \dots, h_n$ are the column
heights of the Ferrers board outlined by $w$.  Then
$$
w = \sum_{k=0}^{n} r_k(q) U^{m-k}D^{n-k},
$$
where
$$
r_k(q) = \df{1}{[n-k]_q\,!} \Delta_q^{n-k}\,P(x)\Big|_{x=0}\;.
$$
\end{corollary}

We define the $q$-Weyl binomial coefficient ${n \brack m}_k$ as the
coefficient of the term $U^{n-m-k}D^{m-k}$ in the product $(D+U)^n$, with the relation
$DU=qUD+1$.

\begin{theorem}\label{T:qWeylBinCoef}
Let $k \le m$ be an integer.  Then  
\begin{align*}
\sum_{B \subseteq [m]\times[n-m]} R_k(B,q) 
  &=  {n\brack m}_k  \\
  &= {n-2k \brack m-k} \;\cdot\; 
     \{x^n y^k\}\; \dfrac{1}
                       {1-x-\dfrac{x^2y}
                               {1-qx-\dfrac{[2]_q\,x^2y}
                                        {1-q^2x-\dfrac{[3]_q\,x^2y}{\cdots}}}} \;,
\end{align*}
where $\{x^n y^k\}\,F(x,y)$ is the coefficient of $x^n y^k$ in the power series $F(x,y)$.
\end{theorem}

\begin{proof}

Each Ferrers board in $[m]\times[n-m]$ is outlined by a unique word with $n-m$ $U$'s and
$m$ $D$'s, and vice versa, so 
$\Sum_{B \subseteq [m]\times[n-m]} R_k(B,q)
  = \displaystyle {n\brack m}_k  $  follows from the definition of the $q$-Weyl
binomial coefficients and Theorem~\ref{T:qmain}.  It remains to explain how the
$q$-Weyl binomial coefficient factors into a $q$-binomial coefficient and a coefficient of
the above continued fraction.

A placement $C$ of $k$ rooks on the Ferrers board $B_w$ corresponds to marking $k$ 
distinct $(D,U)$ pairs on $w$, $D$ to the left of $U$.  This in turn corresponds to
the pair $(w', \pi)$, where $w'$ is the subword of $w$ consisting of the unpaired $n-m-k$
$U$'s and $m-k$ $D$'s, and $\pi$ is the set of $k$ distinct pairs of elements in $[n]$
which keep track of the placement of the marked $(D,U)$ pairs in $w$.

We claim that $\inv(C) = \wt_1(\pi)\cdot \wt_2(w')$ for appropriately chosen weight
functions $\wt_1$ depending only on $\pi$, $\wt_2$ depending only on $w'$.

In terms of the word $w$ with $k$ marked $(D,U)$ pairs,
the statistic $\inv(C)$ is the factor obtained by transforming $w$ into the normally
ordered word $U^{n-m-k}D^{m-k}$, successively replacing $DU$ by $UD$ with a factor of $q$,
deleting each marked pair when the $D$ and the $U$ are adjacent.  We organize the
transformation in the following way.  

First, we bring together and delete the marked
pairs, without disturbing the relative order of the unmarked letters.  During this stage
of transformation, an unpaired letter must commute  with one member of each pair
that it separates.  Also, if two pairs cross---each pair is separated by a member of the
other---then the inner members must commute once.  
Therefore we define $\wt_1(\pi) := q^{\cross(\pi)+\sep(\pi)}$, where
$\cross(\pi)$ is the number of crossings of pairs, and
$\sep(\pi)$ is the number of pairs each unpaired element separates.

Now, being left with the subword $w'$, we commute the letters to their normally
ordered places.  Therefore we define $\wt_2(w') :=q^{|B_{w'}|}$.

Summing over all pairs $(w',\pi)$, we get
\begin{align*}
\displaystyle {n\brack m}_k 
&= \Sum_{w'} q^{|B_{w'}|} \;\cdot\; \Sum_{\pi} q^{\cross(\pi)+\sep(\pi)} \\
&= \displaystyle {n-2k \brack m-k} \;\cdot\; \Sum_{\pi} q^{\cross(\pi)+\sep(\pi)}.
\end{align*}

To find $\sum_{\pi} q^{\cross(\pi)+\sep(\pi)}$, we use a construction of pairings
by weighted Motzkin paths: paths on non-negative integers composed of up steps $i\to i+1$,
constant steps $i\to i$, and down steps $i\to i-1$, weighed respectively by $u_i$,
$c_i$, and $d_i$.   The construction is similar to those in~\cite{F}. A Motzkin path composed of
$n$ steps,
$k$ of which are up steps, constructs a set of $k$ pairs from elements in $[n]$ as follows.

\vs3
\begin{center}
\begin{tabular}{c|p{2in}|c}
   $k$th step   & Construction                       & Weight factor \\ 
   \hline \hline
   $i\to i$     & $k$ is unpaired                    & $c_i=q^i$     \\
   $i\to i+1$   & begin a new pairing with $k$       & $u_i=y$     \\
   $i\to i-1$   & finish a pairing at $k$            & $d_i=[i]_q$    \\ 
\end{tabular}
\end{center}
\vs1

While a path may construct many different pairings, each pairing is constructed by exactly
one Motzkin path. 
With these weights, the weight of a Motzkin path---the product of the weights of
its steps---is the sum $\sum_{\pi} q^{\cross(\pi)+\sep(\pi)}$ over all pairings it
constructs.
 
Since the generating function for weighted Motzkin paths~\cite{F} is the Jacobi continued
fraction
$$
M(x) = \dfrac{1}
             {1-c_0x-\dfrac{u_0d_1x^2}
                           {1-c_1x-\dfrac{u_1d_2x^2}
                                         {1-c_2x-\dfrac{u_2d_3x^2}{\cdots}}}} \;, 
$$
the sum
$\sum_{\pi} q^{\cross(\pi)+\sep(\pi)}$ is the coefficient of $x^n y^k$ in
$$
M(x,y)=\dfrac{1}
             {1-x-\dfrac{[1]_q\,\cdot x^2y}
                        {1-qx-\dfrac{[2]_q\,\cdot x^2y}
                                   {1-q^2x-\dfrac{[3]_q\,\cdot x^2y}{\cdots}}}} \; .
$$

\end{proof}

\section{Further generalizations}\label{S:i-analogue}

We can further generalize the interpretation of normal order coefficients in terms of rook
numbers on Ferrers boards to the algebra with two generators $D$ and $U$ and the commutation
relation $DU-UD=cU^i$, where $i$ is a positive integer.  As in the case of the $q$-Weyl
algebra, this requires a generalization of the rook numbers.

An example of such an algebra is one generated by
$D = \frac{d}{dx}$ and $U$ acting by multiplication by $x^{i+1}$.  It is
easy to check that $DU - UD = (i+1)U^i$.  More generally, if $D = \frac{d}{dx}$ and
$U$ acts as multiplication by a function $U(x)$, then the commutation relation is
$DU-UD=\frac{dU}{dx}$.  Examples of $U$ that give the relation $DU-UD=cU^i$ include $U=e^x$
for $i=1$ and $c=1$, 
$U=\frac{1}{1-x}$ for $i=2$ and $c=1$, and more generally $U=(1-x)^{-1/(i-1)}$ 
for $i\ne 1$ and $c=1/(i-1)$.

Let us consider a word $w$, and find its normal order coefficients. 
Algebraically, we are replacing $DU$ by $UD+cU^i$ and expanding, until all $U$'s are to
the left of all $D$'s in each term.  If we consider the terms as words in the letters $D$ and $U$,
we are either replacing a $DU$ by $UD$, or replacing it with $U^i$, the latter choice weighted by
$c$.  If we consider these choices in terms of the Ferrers board $B_w$ outlined by $w$, then each
time we choose to place a rook, we create $i$ new rows to the right of the square.  This describes
exactly the $i$-row creation rule of rook placement presented by Goldman and Haglund~\cite{GH}.

\begin{example}
We illustrate the $1$-row creation rule by placing three rooks on the Ferrers board $B$ 
with column heights $(3,3,1)$.

\begin{figure}[h]
\begin{center}
\includegraphics[width=0.9in]{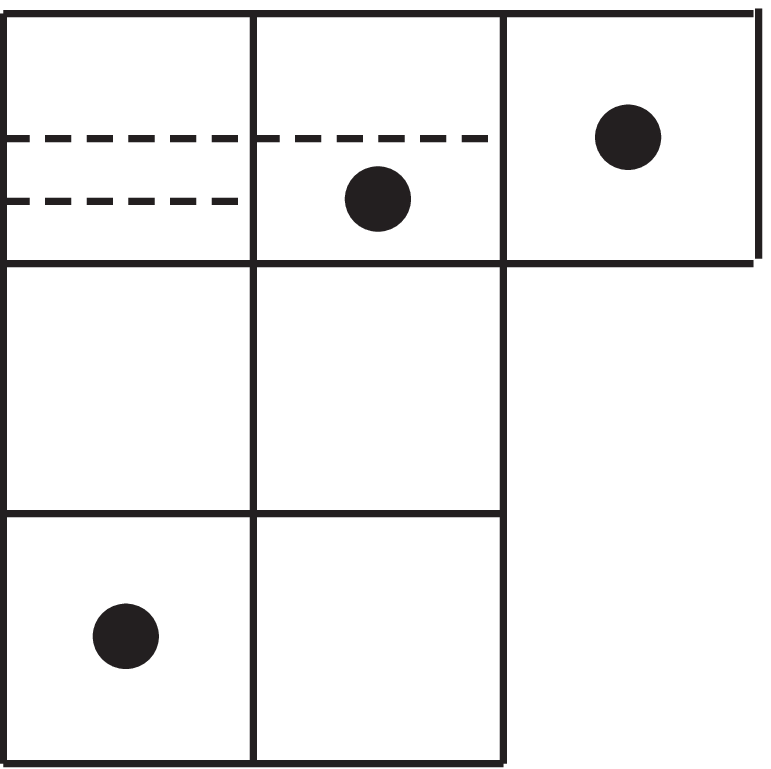}
\end{center}
\end{figure}

We place three rooks on $B$, going from right to left.  Each time we place a rook, we create
a new row to the left of it by splitting the existing row in half.  We consider the bottom
half as the new row, and the top half as the row belonging to the rook just placed.

\end{example}

\begin{definition}
Let $B$ be a Ferrers board.
The $i$-rook number $r_k^{(i)}(B)$ is the number of ways to place $k$ rooks on the
board $B$ going from right to left, creating $i$ new rows to the right of each rook. 
\end{definition}

\begin{theorem}
Let $w$ be an element in the algebra generated by $D$, $U$, with the relation $DU-UD=cU^i$.
If $w$ is a word composed of $n$ $D$'s and $m$ $U$'s, then
$$
w = \sum_{k=0}^{n} c^k\,r_k^{(i)}(B_w) U^{m-k}D^{n-k}.
$$
\end{theorem}

The proof of Theorem~\ref{T:main} easily adapts to this theorem.  Note that since the
placement of rooks on a Ferrers board with $i$-row creation rule always happens from right to left,
it is fortunate that in the proof we consider the rightmost inner square when we decide whether to
place a rook on it.




\bibliographystyle{amsalpha}

\end{document}